\documentclass{article}
\usepackage{amssymb}
\usepackage{amsthm}
\usepackage{amsmath}
\usepackage{hyperref}
\usepackage{listings}
\usepackage{graphicx}
\usepackage{float}
\usepackage{subcaption}
\usepackage{hyperref}\usepackage{chngcntr} 
\graphicspath{{./}}
\usepackage{tikz}
\usetikzlibrary{automata, positioning}
\begin{document}
	\title{A queueing model with servers disguised as customers}
	\author{Prabhu, Janhavi \\
		\texttt{janhaviprabhu27@gmail.com}
		\and
		Hlynka, Myron \\
		\texttt{hlynka@uwindsor.ca}}
	\maketitle
\begin{abstract}
	In this paper, we propose a new queueing model, motivated by the phenomenon of pseudoprogression in cancer, in which the length of a queue appears to increase initially, before reducing to a steady state. We assume that servers arrive to the queue alongside the customers, i.e. `disguised' as customers. We derive the general equations for this model using matrix analytic methods, and demonstrate its behaviour with numerical simulations. 
\end{abstract}

\section{Introduction}
Numerous variations of queueing models have been proposed which take into account heterogenous servers, multiple types of customers and various priorityies. In this paper, we propose a new queueing model, in which servers arrive disguised as customers. The motivation behind this model arose from the documentary
``Jim Allison: Breakthrough'' (\cite{jim}), which spoke about the phenomenon of pseudoprogression in cancer immunotherapy treatments. Pseudoprogression is defined as the increase in tumor size or appearance in new lesions, followed by a decrease in tumor burden. A thorough review of pseudoprogression and its comparison with true progression of a tumor is given in \cite{jia} and \cite{ma} respectively. In a similar manner, our proposed model appears to increase in size initially due to the arrival of servers among the customers, until the servers complete some services, after which the queue size rapidly decreases.

In section (\ref{sec_2}) we review the existing literature on queues with heterogenous customers and servers, as well as matrix geometric solutions for queueing theory. In section (\ref{sec_3}), we introduce notation and the phase diagram for the model, and derive the equations for our model. Section (\ref{sec_4}) involves numerical simulations for the simplest case of our model and considers queue measures for different values of the parameters. Section \ref{sec_5}) gives conclusions. 

\section{Literature Review}\label{sec_2}
There is considerable literature on queueing models with heterogeneous customers. \cite{klim2} presents a queueing model with two types of customers, where the priority changes after a random amount of time. A variation of this appears in  \cite{klim}, which details a multi-server priority queue with marked Markovian heteregenous arrivals. A description of a queueing system with heteregenous customers is given by Collings \cite{coll}, in which he derives the flow through states by performing "cuts" at different points in the system. Customers belong to different groups, having different arrival and service rates, and it is shown that approximating these rates by a single exponential distribution underestimates the queue length. \\
\cite{gurt} presents a queueing system with servers that arrive when the number of customers exceeds a threshold. \cite{yech} put forward a queue with both heteregeneous arrivals and service times. In this model, arrival intensities can "jump" between two different levels. The time spent in each level is exponentially distributed, and all arrivals are assumed to have the same rates during this level. \cite{kim} examines the problem of assigning customers to different servers (having varied service rates) so as to minimize the waiting time of a customer in the case of priority queueing. \cite{li} puts forth a model in which the number of servers in the queue is dynamic, a concept that we use in our model, and depends upon the number of customers in the system. Servers are added or removed on the basis of a search and release process respectively, provided a set of constraints for server utilization are met. \\
We make use of matrix analytic methods for the purpose of representing and solving our model, a detailed review of which is given by \cite{neut}. The procedure to be followed for solving such systems is given in \cite{nels} and \cite{stew}.

\section{Model}\label{sec_3}
\subsection{Notation}

Before presenting our model, we introduce some preliminary notation.
\paragraph*{}
Customers arrive into the system according to a homogeneous Poisson process with rate $\lambda_{c}$, and are served on a first-come-first-serve basis. Service times are independent, exponential random variables with rate $\mu$. In our model, servers arrive in the system along with customers, with rate $\lambda_{sn}$, and leave the system with rate $\mu_{sn}$, where we assume that the arrival and departure of servers can be depend on the number $n$ of customers present in the system. 
\paragraph*{}
We denote the states of the system by $(n,k)$, where  $n\in \left\{0,1,2,3,...\right\}$ represents the number of customers and  $k\in \left\{0,1,2,...,K\right\}$, the number of servers. We restrict the number of servers in the system to $K$, and for the sake of simplicity, we derive the equations for $K = 3$ servers. For each value of $n$, we have $0,1,2$ or $3$ servers in the system. 
\subsection{Phase Diagram}
For each state $(n,k)$, $n$ is the level, while $k$ is the interlevel. For each level, we have $4$ interlevel states. Below is a phase diagram of the process. Row 1 has no servers. Row 2 has 1 server. And so on.  

\makebox[\textwidth][r]{}
		\begin{tikzpicture}[align=center,node distance=1.1cm]
	
	\node[state] (00) {0};
		\node[state, below=of 00] (01) {0'};
	\node[state, right=of 00] (10) {1};
		\node[state, below=of 10] (11) {1'};
		\node[state, below=of 11] (12) {1''};
	\node[state, right=of 10] (20) {2};
		\node[state, below=of 20] (21) {2'};
		\node[state, below=of 21] (22) {2''};
		\node[state, below=of 22] (23) {2'''};
	\node[state, right=of 20] (30) {3};
		\node[state, below=of 30] (31) {3'};
		\node[state, below=of 31] (32) {3''};
		\node[state, below=of 32] (33) {3'''};
	\node[draw = none, right=of 30] (3-n) {$\cdots$};
		\node[draw = none, right=of 31] (3-n1) {$\cdots$};
		\node[draw = none, right=of 32] (3-n2) {$\cdots$};
		\node[draw = none, right=of 33] (3-n3) {$\cdots$};
	\node[state, right=of 3-n] (n0) {n};
		\node[state, below=of n0] (n1) {n'};
		\node[state, below=of n1] (n2) {n''};
		\node[state, below=of n2] (n3) {n'''};
	\node[draw = none, right= of n0] (n+10) {$\cdots$};
		\node[draw = none, right=of n1] (n+11) {$\cdots$};
		\node[draw = none, right=of n2] (n+12) {$\cdots$};
		\node[draw = none, right=of n3] (n+13) {$\cdots$};
	\draw[every loop]
		(00) edge[bend left, auto = left] node {$\lambda_{c}$} (10)
			(01) edge[bend left, auto = left] node {\small $\mu_{s0}$} (00)
			(01) edge[bend left, auto = left] node {$\lambda_{c}$} (11)
		(10) edge[bend left, auto = left] node {$\lambda_{c}$} (20)
			(10) edge[bend left, auto = left] node {\small $\lambda_{s1}$} (11)
			(11) edge[bend left, auto = left] node {\small $\mu_{s1}$} (10)
			(11) edge[bend left, auto = left] node {$\lambda_{c}$} (21)
			(11) edge[bend left, auto = left] node {$\mu$} (01)
			(12) edge[bend left, auto = left] node {\small $\mu_{s1}$} (11)
			(12) edge[bend left, auto = left] node {$\lambda_{c}$} (22)
		(20) edge[bend left, auto = left] node {$\lambda_{c}$} (30)
			(20) edge[bend left, auto = left] node {\small $\lambda_{s2}$} (21)
			(21) edge[bend left, auto = left] node {\small $\mu_{s2}$} (20)
			(21) edge[bend left, auto = left] node {\small $\lambda_{s2}$} (22)
			(21) edge[bend left, auto = left] node {$\lambda_{c}$} (31)
			(21) edge[bend left, auto = left] node {$\mu$} (11)
			(22) edge[bend left, auto = left] node {\small $\mu_{s2}$} (21)
			(22) edge[bend left, auto = left] node {$\lambda_{c}$} (32)
			(22) edge[bend left, auto = left] node {$2\mu$} (12)
			(23) edge[bend left, auto = left] node {\small $\mu_{s2}$} (22)
			(23) edge[bend left, auto = left] node {$\lambda_{c}$} (33)
		(30) edge[bend left, auto = left] node {$\lambda_{c}$} (3-n)
			(30) edge[bend left, auto = left] node {\small $\lambda_{s3}$} (31)
			(31) edge[bend left, auto = left] node {\small $\mu_{s3}$} (30)
			(31) edge[bend left, auto = left] node {\small $\lambda_{s3}$} (32)
			(31) edge[bend left, auto = left] node {$\lambda_{c}$} (3-n1)
			(31) edge[bend left, auto = left] node {$\mu$} (21)
			(32) edge[bend left, auto = left] node {\small $\mu_{s3}$} (31)
			(32) edge[bend left, auto = left] node {\small $\lambda_{s3}$} (33)
			(32) edge[bend left, auto = left] node {$\lambda_{c}$} (3-n2)
			(32) edge[bend left, auto = left] node {$2\mu$} (22)
			(33) edge[bend left, auto = left] node {\small $\mu_{s3}$} (32)
			(33) edge[bend left, auto = left] node {$\lambda_{c}$} (3-n3)
			(33) edge[bend left, auto = left] node {$3\mu$} (23)
		(3-n) edge[bend left, auto = left] node {$\lambda_{c}$} (n0)
			(3-n1) edge[bend left, auto = left] node {$\mu$} (31)
			(3-n1) edge[bend left, auto = left] node {$\lambda_{c}$} (n1)
			(3-n2) edge[bend left, auto = left] node {$2\mu$} (32)
			(3-n2) edge[bend left, auto = left] node {$\lambda_{c}$} (n2)
			(3-n3) edge[bend left, auto = left] node {$3\mu$} (33)
			(3-n3) edge[bend left, auto = left] node {$\lambda_{c}$} (n3)
			(n0) edge[bend left, auto = left] node {$\lambda_{c}$} (n+10)
			(n0) edge[bend left, auto = left] node {\small $\lambda_{sn}$} (n1)
			(n1) edge[bend left, auto = left] node {\small $\mu_{sn}$} (n0)
			(n1) edge[bend left, auto = left] node {\small $\lambda_{sn}$} (n2)
			(n1) edge[bend left, auto = left] node {$\lambda_{c}$} (n+11)
			(n1) edge[bend left, auto = left] node {$\mu$} (3-n1)
			(n2) edge[bend left, auto = left] node {\small $\mu_{sn}$} (n1)
			(n2) edge[bend left, auto = left] node {\small $\lambda_{sn}$} (n3)
			(n2) edge[bend left, auto = left] node {$\lambda_{c}$} (n+12)
			(n2) edge[bend left, auto = left] node {$2\mu$} (3-n2)
			(n3) edge[bend left, auto = left] node {\small $\mu_{sn}$} (n2)
			(n3) edge[bend left, auto = left] node {$\lambda_{c}$} (n+13)
			(n3) edge[bend left, auto = left] node {$3\mu$} (3-n3)
			(n+11) edge[bend left, auto = left] node {$\mu$} (n1)
			(n+12) edge[bend left, auto = left] node {$2\mu$} (n2)
			(n+13) edge[bend left, auto = left] node {$3\mu$} (n3);
		
\end{tikzpicture}

The dashes beside the numbers in the above diagram represent the number of servers. Assume that no new servers can arrive if the number of servers is greater than or equal to the number of current customers. Assume that if the number of servers is one more than the number of customers, the remaining servers are dormant until the extra server leaves, or another customer arrives. This can be seen in the transition $(1,2)\rightarrow(1,1)$.  We will now derive the time dependent probability equations for the model.

\subsection{Matrix Analytic Solution of the Model}
We now proceed to derive the equations for the most  general case of the model. Order of the states is $0,0',1,1',1'',2,2',2'',2''',\dots$. The generator matrix for the system is given by:
\paragraph*{}
	$Q = 
	\begin{pmatrix}
		-\lambda_{c} & \lambda_{c} & 0 & 0 & 0 & 0 & 0 & \dots\\
		\mu_{s0} & - \mu_{s0} - \lambda_{c}  & 0 & \lambda_{c} & 0 & 0 & 0 & \dots \\
		0 & 0 & -\lambda_{s1} - \lambda_{c} & \lambda_{s1} & 0 & \lambda_{c} & 0 & \dots \\
		0 & \mu & \mu_{s1} & -\mu_{s1} - \mu - \lambda_{c} & 0 & 0 & \lambda_{c} & \dots\\
		0 & 0 & 0 & \mu_{s1} & -\mu_{s1} - \lambda_{c} & 0 & 0 & \dots \\
		\vdots & \vdots & \vdots & \vdots & \vdots & \vdots & \vdots & \ddots
	\end{pmatrix}$ 
\paragraph*{}
We can write Q in the form of a block matrix as follows: 
\begin{equation}
	Q = \begin{pmatrix}
		B_{00} & B_{01} & 0 & 0 & 0 & 0 & 0 & \dots \\
		B_{10} & B_{11} & B_{12} & 0 & 0 & 0 & 0 & \dots \\
		0 & B_{21} & B_{22} & A_{1} & 0 & 0 & 0 & \dots \\
		0 & 0 & A_{0} & A_{3} & A_{1} & 0 & 0 & \dots \\
		0 & 0 & 0 & A_{0} & A_{4} & A_{1} & 0 & \dots \\
		0 & 0 & 0 & 0 & A_{0} & A_{5} & A_{1} & \dots \\
		\vdots & \vdots & \vdots & \vdots & \vdots & \vdots & \vdots & \ddots \\
	\end{pmatrix}
\end{equation}
where:
\begin{equation}
	\label{eqn_2}
	\begin{aligned}
		B_{00} &= \begin{bmatrix}
			-\lambda_{c} & 0 \\
			\mu_{s0} & -\mu_{s0} - \lambda_{c}
		\end{bmatrix}, \qquad
		B_{01} = \begin{bmatrix}
			\lambda_{c} & 0 & 0 \\
			0 & \lambda_{c} & 0
		\end{bmatrix}, \\[1em]
		B_{10} &= \begin{bmatrix}
			0 & 0 \\
			0 & \mu \\
			0 & 0
		\end{bmatrix}, \qquad
		B_{11} = \begin{bmatrix}
			-\lambda_{s1}-\lambda_c & \lambda_{s1} & 0 \\
			\mu_{s1} & -\mu - \mu_{s1} - \lambda_{c}  & 0 \\
			0 & \mu_{s1} & -\mu_{s1} - \lambda_{c}
		\end{bmatrix}, \\[1em]
		B_{12} &= \begin{bmatrix}
			\lambda_{c} & 0 & 0 & 0 \\
			0 & \lambda_{c} & 0 & 0 \\
			0 & 0 & \lambda_{c} & 0
		\end{bmatrix}, \qquad
		B_{21} = \begin{bmatrix}
			0 & 0 & 0 \\
			0 & \mu & 0 \\
			0 & 0 & 2\mu \\
			0 & 0 & 0
		\end{bmatrix} \qquad \text{and}\\[1em]
		B_{22} &= \begin{bmatrix}
			-\lambda_{s2}-\lambda_c & \lambda_{s2} & 0 & 0 \\
			\mu_{s2} & -\mu_{s2} - \mu - \hat{\lambda} & \lambda_{s2} & 0 \\
			0 & \mu_{s2} & -\mu_{s2} - 2\mu - \lambda_{s2} - \lambda_{c} & 0 \\
			0 & 0 & \mu_{s2} & -\mu_{s2} - \lambda_{c} 
		\end{bmatrix}
	\end{aligned}
\end{equation}
are known as the boundary conditions. We also have the matrices:
\begin{equation}
	\label{eqn_3}
	\begin{aligned}
		A_{1} &= \begin{bmatrix}
			\lambda_{c} & 0 & 0 & 0 \\
			0 & \lambda_{c} & 0 & 0 \\
			0 & 0 & \lambda_{c} & 0 \\
			0 & 0 & 0 & \lambda_{c}
		\end{bmatrix} \\[1em]
		A_{0} &= \begin{bmatrix}
			0 & 0 & 0 & 0 \\
			0 & \mu & 0 & 0 \\
			0 & 0 & 2\mu & 0 \\
			0 & 0 & 0 & 3\mu
		\end{bmatrix}
	\end{aligned}
\end{equation}
and: 
\begin{equation*}
	A_{n} = \begin{bmatrix}
		-\lambda_{sn} - \lambda_{c} & \lambda_{sn} & 0 & 0 \\
		\mu_{sn} & -\mu - \mu_{sn} - \lambda_{sn} - \lambda_{c} & \lambda_{sn} & 0 \\
		0 & \mu_{sn} & -2\mu - \mu_{sn} - \lambda_{sn} - \lambda_{c} & \lambda_{sn} \\
		0 & 0 & \mu_{sn} & -3\mu - \mu_{sn} - \lambda_{c}
	\end{bmatrix}
\end{equation*} \\
where $n \geq 3$. Notice that the matrix $A_{0}$ represents the movement to a lower level and $A_{1}$ represents the movement to a higher level.

Let $\pi = (\pi_{0}, \pi_{1}, \pi_{2}, ...)$ denote the limiting distribution of the system, where $\pi_{i}$ denotes the limiting distribution of the $i^{th}$ level. Then, we have the equations:
\begin{align}
	\label{eqn_4}
	\pi Q &= 0 \qquad \text{and} \\
	\label{eqn_5}
	\pi \textbf{e} &= 1
\end{align}
where \textbf{e} denotes a vector of 1's. Solving the above two equations gives us a set of boundary conditions:
\begin{equation} \label{eqn_6}
	\begin{aligned}
		\pi_{0}B_{00} + \pi_{1}B_{10} = 0 \\
		\pi_{0}A_{1} + \pi_{1}B_{11} + \pi_{2}B_{21} = 0 \\
		\pi_{1}A_{1} + \pi_{2}B_{22} + \pi_{3}A_{0} = 0
	\end{aligned}
\end{equation}
For $i \geq 3$, we have the following recursive relation:
\begin{equation} \label{eqn_7}
	\pi_{i-1}A_{1} + \pi_{i}A_{i} + \pi_{i+1}A_{0} = 0
\end{equation}
To solve the above equations, we need to express $A_{n}$ in terms of the values we know. In the following section, we take the simplest case, that is, when $\lambda_{sn} = \lambda_{s}$ and $\mu_{sn} = \mu_{s}$ and derive the properties of the corresponding queue.
\subsection{Simplest Case: \texorpdfstring{$\lambda_{sn} = \lambda_{s}$ and $\mu_{sn} = \mu_{s}$}{Constant parameters}}
In this case, $A_{n}$ becomes a constant matrix given by:
\begin{equation}
	\label{eqn_8}
	A_{2} = \begin{bmatrix}
		-\hat{\lambda} & \lambda_{s} & 0 & 0 \\
		\mu_{s} & -\mu - \mu_{s} - \hat{\lambda} & \lambda_{s} & 0 \\
		0 & \mu_{s} & -2\mu - \mu_{s} - \hat{\lambda} & \lambda_{s} \\
		0 & 0 & \mu_{s} & -3\mu - \mu_{s} - \lambda{c} 
	\end{bmatrix}
\end{equation}
where $\hat{\lambda} = \lambda_{c} + \lambda_{s}$. Our generator matrix becomes:
\begin{equation}
	\label{eqn_9}
	Q = \begin{pmatrix}
		B_{00} & B_{01} & 0 & 0 & 0 & 0 & ) & \dots \\
		B_{10} & B_{11} & B_{12} & 0 & 0 & 0 & 0 & \dots \\
		0 & B_{21} & B_{22} & A_{1} & 0 & 0 & 0 & \dots \\
		0 & 0 & A_{0} & A_{2} & A_{1} & 0 & 0 & \dots \\
		0 & 0 & 0 & A_{0} & A_{2} & A_{1} & 0 & \dots \\
		0 & 0 & 0 & 0 & A_{0} & A_{2} & A_{1} & \dots \\
		\vdots & \vdots & \vdots & \vdots & \vdots & \vdots & \vdots & \ddots \\
	\end{pmatrix}
\end{equation}
The recurrence relation equation (\ref{eqn_7}) reduces to the form:
\begin{equation}
	\label{eqn_10}
	\pi_{i-1}A_{1} + \pi_{i}A_{2} + \pi_{i+1}A_{0} = 0
\end{equation}
Notice that in the matrix Q, once we reach the third level $\left(n \geq 3\right)$, the matrix starts repeating, precisely like a Quasi Birth-and-Death process. We can now solve for the limiting distributions. 
\subsubsection{Condition for stability}
Let the generator matrix A be given by $A = A_{0} + A_{1} + A_{2}$. This forms the generator for the repeating section of the system. To derive a condition for the stability of the queue, we first find the stationary distribution for A, given by $\pi_{A} = (\pi_{A0}, \pi_{A1}, \pi_{A2}, \pi_{A3})$. For this, we use the equations: 
\begin{equation*}
	\begin{aligned}
		\pi_{A}A &= 0 \qquad \text{and} \\
		\sum_{i = 0}^{3}\pi_{Ai} &= 1
	\end{aligned}
\end{equation*}
Then, we have:
\begin{equation}
	\label{eqn_11}
	\begin{aligned}
		(\pi_{A0}, \pi_{A1}, \pi_{A2}, \pi_{A3}) \cdot 
		\begin{pmatrix}
			-\lambda_{s} & \lambda_{s} & 0 & 0 \\
			\mu_{s} & -\mu_{s} - \lambda_{s} & \lambda_{s} & 0 \\
			0 & \mu_{s} & -\mu_{s} - \lambda_{s} & \lambda_{s} \\
			0 & 0 & \mu_{s} & -\mu_{s}
		\end{pmatrix} &= 0 \\
	\end{aligned}
\end{equation}
Let $\rho_s=\lambda_s/\mu_s$. Solving for $(\pi_{A0}, \pi_{A1}, \pi_{A2}, \pi_{A3})$,
\begin{equation}
	\label{eqn_12}
	\begin{aligned}
		-\lambda_{s}\pi_{A0} + \mu_{s}\pi_{A1} &= 0 \\
		\Rightarrow \pi_{A1} &= \rho_{s}\pi_{A0} \\
		\lambda_{s} - \left(\mu_{s} + \lambda_{s}\right)\pi_{A1} + \mu_{s}\pi_{A2} &= 0 \\
		\Rightarrow \pi_{A2} = \rho_{s}\pi_{A1} &= \rho_{s}^{2}\pi_{A0} \\
		\lambda_{s}\pi_{A1} - \left(\mu_{s} + \lambda_{s}\right)\pi_{A2} + \mu_{s}\pi_{A3} &= 0 \\
		\Rightarrow \pi_{A3} = \rho_{s}\pi_{A2} &= \rho_{s}^{3}\pi_{A0}\\	
	\end{aligned}
\end{equation}
From Equation (\ref{eqn_12}), and the normalization equation, we get:
\begin{equation}
	\label{eqn_13}
	\sum_{i=0}^{3}\pi_{Ai} = 1 \Rightarrow \pi_{A0} = \dfrac{1-\rho_{s}}{1-\rho_{s}^{4}}, \lambda_{s} \neq \mu_{s}
\end{equation}
For ergodicity, we require the \textit{drift} to lower levels to be greater than that to higher levels. That is, the system will be more inclined to go to lower level.
\begin{equation}
	\label{eqn_14}
	\begin{aligned}
		\pi_{A}A_{1}\textbf{e} &< \pi_{A}A_{0}\textbf{e} \\[1em]
		\pi_{A}\cdot
		\begin{bmatrix}
			\lambda_{c} & 0 & 0 & 0 \\
			0 & \lambda_{c} & 0 & 0 \\
			0 & 0 & \lambda_{c} & 0 \\
			0 & 0 & 0 & \lambda_{c}
		\end{bmatrix}
		\begin{bmatrix}
			1 \\ 1 \\ 1 \\ 1
		\end{bmatrix} &< \pi_{A}\cdot
		\begin{bmatrix}
			0 & 0 & 0 & 0 \\
			0 & \mu & 0 & 0 \\
			0 & 0 & 2\mu & 0 \\
			0 & 0 & 0 & 3\mu
		\end{bmatrix}
		\begin{bmatrix}
			1 \\ 1 \\ 1 \\ 1
		\end{bmatrix}\\[1em]
		\Rightarrow \lambda_{c}\sum_{i=0}^{3}\pi_{Ai} &< \left(0, \mu\pi_{A1}, 2\mu\pi_{A2}, 3\mu\pi_{A3}\right) \begin{bmatrix}
		1 \\ 1 \\ 1 \\ 1
		\end{bmatrix}\\
	\end{aligned}
\end{equation}
This gives us our ergodicity condition:
\begin{equation}
	\label{eqn_15}
	\rho_{c} < \dfrac{\rho_{s}(1 + 2\rho_{s} + 3\rho_{s}^{2})}{(1 + \rho_{s})(1 + \rho_{s}^{2})} \quad \text{where } \rho_{c} = \lambda_{c}/\mu
\end{equation}
\subsubsection{Solving for Limiting distribution}
	Once we check whether the queue will be stable for the given values of $\lambda_{c}, \mu, \lambda_{s}$ and $\mu_{s}$, we can find the limiting distribution. From Equation (\ref{eqn_10}), since all three matrices $A_{0}, A_{1}$ and $A_{2}$ are constant, we can assume that for $i \geq 3$, $\pi_{i} = \pi_{i-1}R$, where R is a constant matrix. Repeatedly substituting in the recurrence relation, we get:
	\begin{equation}
		\label{eqn_16}
		\pi_{i} = \pi_{2}R^{i-2}, i \geq 3
	\end{equation}
	The simplest method to solve for R is to iterate numerically until the elements of the matrix differ by less than a specified small value. Substituting the value of $\pi_{i-1}, \pi_{i}$ and $\pi_{i+1}$, the iterative formula is found as follows:
	\begin{equation*}
		\begin{split}
			R^{2}A_{0} + RA_{2} + A_{1} &= 0 \\
			\Rightarrow R^{2}A_{0}A_{2}^{-1} + R + A_{1}A_{2}^{-1} &= 0 \\
			R = -A_{1}A_{2}^{-1} - R^{2}A_{0}A_{2}^{-1} &= -V - R^{2}W \quad \text{where} \\
			R_{(0)} = 0, \quad R_{(k+1)} &= -V - R_{(k)}^{2}W, \enspace k=0,1,2,...
		\end{split}	
	\end{equation*}
	We now solve the boundary equations (\ref{eqn_6}) for $\pi_{0}, \pi_{1}$ and $\pi_{2}$. Substituting the value of $\pi_{3} = \pi_{2}R$, we can represent the equations in matrix form as:
	\begin{equation*}
		\begin{split}
			(\pi_{0}, \pi_{1}, \pi_{2}) \begin{pmatrix}
				B_{00} & B_{01} & 0 \\
				B_{10} & B_{11} & B_{12} \\
				0 & B_{21} & B_{22} + RA_{0}
			\end{pmatrix} = (0, 0 | 0, 0, 0| 0, 0, 0, 0)
		\end{split}
	\end{equation*}
	Note that $\pi_{0} = (\pi_{00}, \pi_{01}), \pi_{1} = (\pi_{10}, \pi_{11}, \pi_{12})$ and $\pi_{2} = (\pi_{20}, \pi_{21}, \pi_{22}, \pi_{23})$. We take $\pi_{00} = 1$, that is, we always start in the $\left(0,0\right)$ state, and solve for the remaining 8 values. We can rewrite the above equation as follows:
	\begin{equation}
		\label{eqn_17}
		\begin{pmatrix}
			B_{00}^\intercal & B_{10}^\intercal & 0 \\[0.5em]
			B_{01} ^\intercal & B_{11}^\intercal & B_{21}^\intercal \\[0.5em]
			0 & B_{12}^\intercal & \left(B_{22} + RA_{0}\right)^\intercal
		\end{pmatrix}
	\begin{pmatrix}
		\pi_{0}^\intercal \\[0.5em] \pi_{1}^\intercal \\[0.5em] \pi_{2}^\intercal
	\end{pmatrix} = \begin{pmatrix}
	0 \\ 0 \\ \hline 0 \\ 0 \\ 0 \\ \hline 0 \\ 0 \\ 0 \\ 0
	\end{pmatrix}
	\end{equation}
Once we have values for $\pi_{0}, \pi_{1}$ and $\pi_{2}$, we normalize them to get probabilities:
\begin{equation}
	\label{eqn_18}
	\begin{aligned}
		\pi \textbf{e} &= 1 \\
		\alpha &= \pi_{0}\textbf{e} + \pi_{1}\textbf{e} + \pi_{2} \textbf{e} + \sum_{i=1}^{\infty}\pi_{2}R^{i} \\
		\Rightarrow \alpha &= \pi_{0}\textbf{e} + \pi_{1}\textbf{e} + \sum_{i=0}^{\infty}\pi_{2}R^{i} \\
		\alpha &= \pi_{0}\textbf{e} + \pi_{1}\textbf{e} + \pi_{2}\left(I-R\right)^{-1}
	\end{aligned}
\end{equation}
In the above equations, we take \textbf{e} merely to be a vector of ones of required dimension for each term. For preciseness, we can take $\mathbf{e_{0}} = [1 \enspace 1]^\intercal$ and $\mathbf{e_{1}}= [1 \enspace 1 \enspace 1]^\intercal$ and write $\pi_{0}\mathbf{e_{0}} + \pi_{1}\mathbf{e_{1}}$ instead of $\pi_{0}\textbf{e} + \pi_{1}\textbf{e}$. \\
$\alpha$ is a scalar, and we divide $\pi_{0}, \pi_{1}$ and $\pi_{2}$ to normalize the probabilities. Once we have the boundary limiting distributions we can find any $\pi_{n} = \pi_{2}R^{n-2}$. In the next section, we find the properties of the queue.

\subsubsection{Properties of the queue}
\paragraph{Expected queue length}:
For our model, the expected queue length includes the number of servers, since an observer of the queue cannot distinguish between a customer and a server. Before we find an expression for the expected queue length, we first define:
\begin{align*}
	\mathbf{v_{0}} &= \begin{bmatrix}
		0 \\ 1
	\end{bmatrix} & 
	\mathbf{v_{1}} &= \begin{bmatrix}
		0 \\ 1 \\ 2
	\end{bmatrix} &
	\mathbf{v} &= \begin{bmatrix}
		0 \\ 1 \\ 2 \\ 3
	\end{bmatrix}
\end{align*}
Then, the expression for expected queue length will be given by:
\begin{equation*}
	\begin{aligned}
		E(L) &= \sum_{n=0}^{\infty}\sum_{k}\left(n+k\right)p_{(n,k)} \\
			 &= \pi_{0}\mathbf{v_{0}} + \pi_{1}(\mathbf{e_{1}} + \mathbf{v_{1}}) + \sum_{n=2}^{\infty}\sum_{k}(n+k)p_{(n,k)} \\
			 &= \pi_{0}\mathbf{v_{0}} + \pi_{1} (\mathbf{e_{1}} + \mathbf{v_{1}}) + \sum_{n=2}^{\infty}\pi_{2}R^{n-2}\textbf{e} + \sum_{n=2}^{\infty}\pi_{2}R^{n-2}\textbf{v} \\
			 &= \pi_{0}\mathbf{v_{0}} + \pi_{1} (\mathbf{e_{1}} + \mathbf{v_{1}}) + \pi_{2}(\sum_{n=2}^{\infty}nR^{n})\textbf{e} + \pi_{2}(\sum_{n=0}^{\infty}R^{n})\textbf{v} \\
	\end{aligned}
\end{equation*}
By manipulating the third term and changing the limits, we can obtain our final equation:
\begin{multline}
	\label{eqn_19}
	E(L) = \pi_{0}\mathbf{v_{0}} + \pi_{1} (\mathbf{e_{1}} + \mathbf{v_{1}}) + \pi_{2}[R(I-R)^{2} + 2(I-R)^{-1}]\textbf{e} + \pi_{2}(I-R)^{-1}\textbf{v}
\end{multline}

\paragraph*{Expected Number of Customers Waiting in Queue} :
We denote this value by $E(L_{n})$. The derivation for the expected number of customers is similar to that for the expected queue length.
\begin{equation*}
	\begin{aligned}
		E(L_{n}) &= \sum_{n=1}\sum_{k}(n-k)p(n,k) \\
				 &= \sum_{n=1}^{3}\sum_{k=0}^{n}(n-k)p(n,k) + \sum_{n=4}^{\infty}\sum_{k=0}^{3}(n-k)p(n,k) \\
				 &= \sum_{n=1}^{3}\sum_{k=0}^{n}(n-k)p(n,k) + \sum_{n=4}^{\infty}\sum_{k=0}^{3}np(n,k) - \sum_{n=4}^{\infty}\sum_{k=0}^{3}kp(n,k) \\
				 &= \sum_{n=1}^{3}\sum_{k=0}^{n}(n-k)p(n,k) + \pi_{2}(\sum_{n=4}^{\infty}nR^{n-2}\textbf{e} - \sum_{n=4}^{\infty}R^{n-2}\textbf{v})
	\end{aligned}
\end{equation*}
By rearranging the terms and changing the limits, we get the following expression:
\begin{multline}
	\label{eqn_20}
	E(L_{n}) = \pi_{10} + (2\pi_{20} + \pi_{21}) + (3\pi_{30} + 2\pi_{31} + \pi_{32}) + \pi_{2}[R(I-R)^{-2} \\ + 2(I-R)^{-1} - (2I + 3R)]\textbf{e} - \pi_{2}[(I-R)^{-1} - (I+R)]\textbf{v}
\end{multline}
The expected system time can  be found by using Little's Law: \\
\begin{center}
	$E(L_{n}) = \lambda_{c}E(W) \Rightarrow E(W) = \dfrac{E(L_{n})}{\lambda_{c}})$
\end{center} 
\paragraph*{Delay Probability}: 
The probability a customer has to wait, known as \textit{delay probability}, is denoted by $\Pi_{w}$, and it is derived as follows:
\begin{equation}
	\label{eqn_21}
	\begin{aligned}
		\Pi_{w} &= \sum_{k}\sum_{n=k+1}^{\infty}p(n,k) \\
				&= \sum_{n=1}^{3}\sum_{k=0}^{n}p(n,k) + \sum_{n=4}^{\infty}\pi_{n}\textbf{e} \\
				&= (\pi_{10} + \pi_{11}) + (\pi_{20} + \pi_{21} + \pi_{22}) + \pi_{2}((I-R)^{-1} - I)\textbf{e}
	\end{aligned}
\end{equation}
\section{Numerical Simulations}\label{sec_4}
We simulate the model for different values of $\lambda_{c}, \lambda_{s}, \mu$ and $\mu_{s}$ and display the results of the model below. \\
\begin{center}
	\begin{tabular} {||c|c|c|c|c|c||}
		\hline
		\multicolumn{6}{|c|}{$\lambda_{c} = 1, \lambda_{s} = 2$}\\ 
		\hline
		$\mu$ & $\mu_{s}$ & E(L) & E($L_{n}$) & E(W) & $\Pi_{w}$ \\
		\hline \hline
		1 & 1 & 2.609 & 0.287 & 2.609 & 0.399 \\
		1 & 2 & 1.533 & 0.426 & 1.533 & 0.375 \\
		1 & 3 & 1.166 & 0.498 & 1.166 & 0.339 \\
		2 & 1 & 1.478 & 0.204 & 1.478 & 0.282 \\
		2 & 2 & 0.894 & 0.270 & 0.894 & 0.250 \\
		2 & 3 & 0.710 & 0.302 & 0.710 & 0.225 \\
		3 & 1 & 1.042 & 0.150 & 1.042 & 0.213 \\
		3 & 2 & 0.669 & 0.193 & 0.669 & 0.186 \\
		3 & 3 & 0.543 & 0.214 & 0.543 & 0.168 \\
		\hline
		\multicolumn{6}{|c|}{$\lambda_{c} = 2, \lambda_{s} = 1$}\\ 
		\hline
		$\mu$ & $\mu_{s}$ & E(L) & E($L_{n}$) & E(W) & $\Pi_{w}$ \\
		\hline \hline
		1 & 1 & 4.474 & 0.756 & 0.378 & 0.557 \\
		1 & 2 & 4.055 & 1.627 & 0.813 & 0.657 \\
		1 & 3 & 4.152 & 2.450 & 1.250 & 0.713 \\
		2 & 1 & 3.559 & 0.764 & 0.382 & 0.490 \\
		2 & 2 & 3.211 & 1.493 & 0.747 & 0.579 \\
		2 & 3 & 3.282 & 2.096 & 1.048 & 0.621 \\
		3 & 1 & 2.934 & 0.732 & 0.366 & 0.440 \\
		3 & 2 & 2.658 & 1.311 & 0.655 & 0.509 \\
		3 & 3 & 2.709 & 1.757 & 0.879 & 0.543 \\	
		\hline	
	\end{tabular}
\end{center}
We can make a few observations from the above values:
\begin{itemize}
	\item As $\mu$ increases, the values of $E(L_{n})$ decrease. That is, as the service rate increases, fewer customers are found in the queue. Naturally, the expected system time $E(W)$ will also decrease if customers are served faster. The delay probability $\Pi_{w}$ also decreases with an increase of $\mu$.
	\item An increase of $\mu_{s}$ results in an increased number of customers in the system ($E(L_{n})$). This is consistent with the assumption that as servers leave the queue at a faster rate, more customers will be left waiting. Note that when $\lambda_{c} < \lambda_{s}$, the delay probability and waiting time decrease even if $\mu_{s}$ increases. This could be attributed to the fact that though the rate of departure of servers increases, their rate of arrival makes up for it. Notice that when $\lambda_{c} > \lambda_{s}$, the expected system time and delay probability increase with $\mu_{s}$.
\end{itemize}
A visualization of the above observations is shown in the plots below. We have kept $\lambda_{s} = 3$ fixed for ease of plotting.
\begin{figure}[htp]
	\includegraphics[scale=0.4]{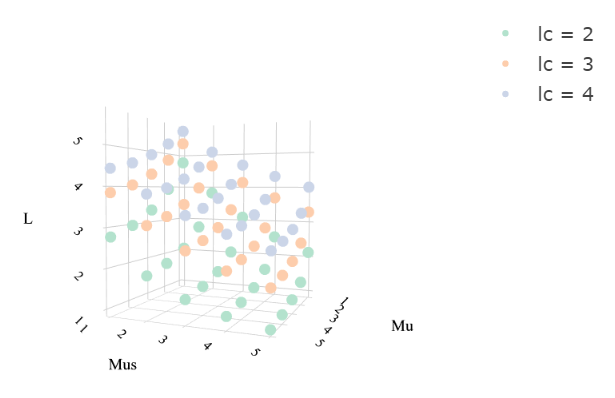}
	\centering
	\caption{Value of E(L) with varying $\mu, \mu_{s}$ and $\lambda_{c}$}
	\label{fig1}
\end{figure}
\begin{figure}[htp]
	\includegraphics[scale = 0.4]{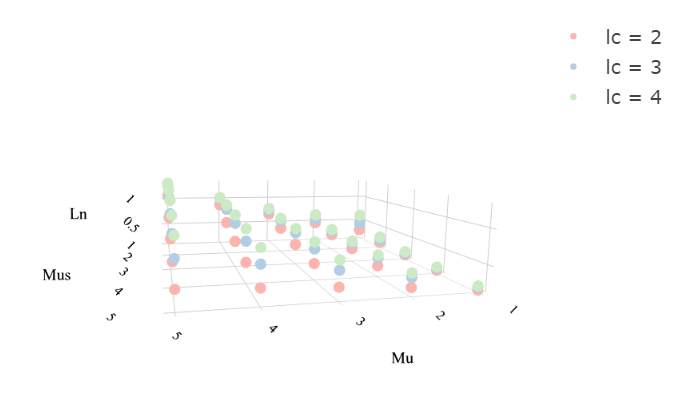}
	\centering
	\caption{Value of E($L_{n}$) with varying $\mu, \mu_{s}$ and $\lambda_{c}$}
	\label{fig2}
\end{figure}
\begin{figure}[htp]
	\includegraphics[scale = 0.4]{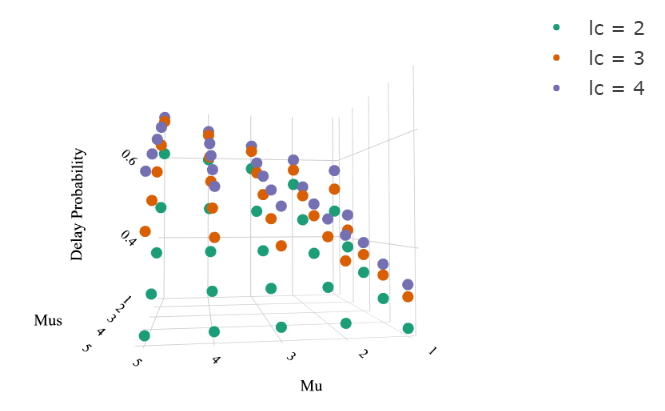}
	\centering
	\caption{Value of $\Pi_{w}$ with varying $\mu, \mu_{s}$ and $\lambda_{c}$}
	\label{fig3}
\end{figure}

\paragraph*{Example 1:} $\lambda_{c} = 1, \mu_ = 2, \lambda_{s} = 3, \mu_{s} = 4$:

Here, $\rho_{c} = 0.5, \rho_{s} = 0.75$ and $(\rho_{s}(1 + 2\rho_{s} + 3\rho_{s}^{2}))/((1 + \rho_{s})(1 + \rho_{s}^{2})) = 1.149 $, hence the ergodicity condition is satisfied - the queue will be stable. 
\begin{align*}
	B_{00} &=
	\begin{bmatrix}
		-1 & 1 \\ 4 & -5
	\end{bmatrix}, &
	B_{01} &= \begin{bmatrix}
		0 & 0 & 0 \\ 0 & 1 & 0
	\end{bmatrix}, &
	B_{10} &= \begin{bmatrix}
		1 & 0 \\ 0 & 2 \\ 0 & 0
	\end{bmatrix} \\[0.75em]
	B_{11} &= \begin{bmatrix}
		-4 & 3 & 0 \\ 4 & -7 & 0 \\ 0 & 4 & -5
	\end{bmatrix}, &
	B_{12} &= \begin{bmatrix}
		1 & 0 & 0 & 0 \\ 0 & 1 & 0 & 0 \\ 0 & 0 & 1 & 0
	\end{bmatrix}, &
	B_{21} &= \begin{bmatrix}
		0 & 0 & 0 \\ 0 & 2 & 0 \\ 0 & 0 & 4 \\ 0 & 0 & 0
	\end{bmatrix} \\[0.75em] 
	B_{22} &= \begin{bmatrix}
		-4 & 3 & 0 & 0 \\ 4 & -10 & 3 & 0 \\ 0 & 4 & -12 & 3 \\ 0 & 0 & 4 & -5
	\end{bmatrix}, & A_{0} &= 
	\begin{bmatrix}
	0 & 0 & 0 & 0 \\ 0 & 2 & 0 & 0 \\ 0 & 0 & 4 & 0 \\ 0 & 0 & 0 & 6
	\end{bmatrix}, & A_{1} &= I \\[0.75em]
	A_{2} &= 
	\begin{bmatrix}
		-4 & 3 & 0 & 0 \\ 4 & -10 & 3 & 0 \\ 0 & 4 & -12 & 3 \\ 0 & 0 & 4 & -11
	\end{bmatrix}
\end{align*}
Using Equation (\ref{eqn_13}), we get $\pi_{A} = (0.366, 0.274, 0.206, 0.154)$ with condition (\ref{eqn_15}) indicating that the system is stable. \\
We found R iteratively, till the difference between elements of successive iterations are less than $10^{-12}$. This gave R as:
\begin{center}
	$R = \begin{bmatrix}
		0.3771186 & 0 & 0 & 0 \\
		0 & 0.1694915 & 0 & 0 \\
		0 & 0 & 0.1087571 & 0 \\
		0 & 0 & 0 & 0.1016949
	\end{bmatrix}$
\end{center}
Now we proceed to solve the boundary conditions (\ref{eqn_6}). We set $\pi_{00} = 1$, that is, we start at the (0,0) state, and find the remaining 8 variables. Given the equation in the form A\textbf{x} = \textbf{0}, we can re-write this as A'\textbf{x'} = \textbf{b}, where $\mathbf{b} = -A[,1]$, and A' is the matrix obtained by taking the submatrix A[2:9,2:9] of A (since we require 8 variables, we need only 8 equations, hence we take an 8x8 submatrix of A). Once we solve the matrix equation, we normalize the probabilities using Equation (\ref{eqn_18}. Our final probabilities are:
\begin{align*}
	\pi_{0} &= (0.6557041, 0.1519753) \\
	\pi_{1} &= (0.052086275, 0.052086275, 0.005271742) \\
	\pi_{2} &= (0.030662973, 0.017641404, 0.006589678, 0.004503371)
\end{align*}
Once we have $\pi_{0}, \pi_{1} \text{ and } \pi_{2}$, we can find all other probabilities of the system as well as the properties: \\
\begin{itemize}
	\item The total expected queue length of the system E(L) = 0.576480174113887 
	\item The expected waiting time of a customer W = 0.576480174113887 
	\item The expected queue length (customers only) E($L_{n}$) = 0.206882719667415
	\item The probability an arriving customer has to wait $\Pi_{w}$ = 0.182545493796219
\end{itemize}
A plot of the simulation of the instantaneous queue length is given below:
\begin{figure}[htp]
	\includegraphics[scale = 0.35]{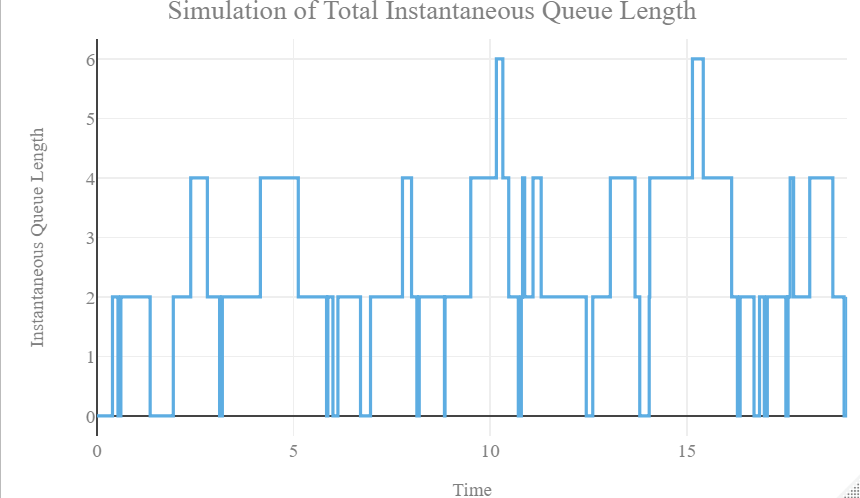}
	\centering
	\caption{Instantaneous total queue length}
	\label{fig4}
\end{figure}
\begin{figure}[htp]
	\includegraphics[scale = 0.35]{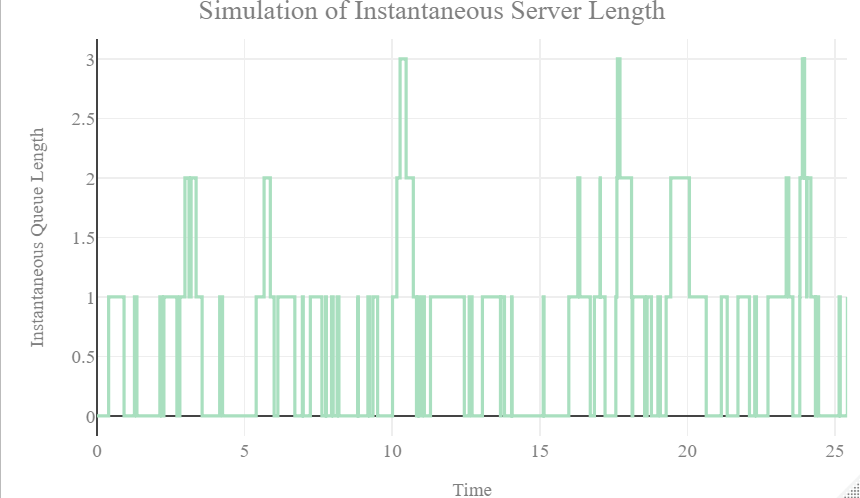}
	\centering
	\caption{Instantaneous server queue length}
	\label{fig5}
\end{figure}
\begin{figure}[htp]
	\includegraphics[scale = 0.35]{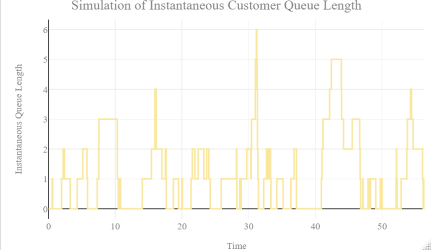}
	\centering
	\caption{Instantaneous customer queue length}
	\label{fig6}
\end{figure}

\paragraph*{Example 2:} $\lambda_{c} = 2, \mu_ = 1, \lambda_{s} = 1, \mu_{s} = 2$:

Here, $\rho_{c} = 2, \rho_{s} = 0.5$ and $\dfrac{\rho_{s}(1 + 2\rho_{s} + 3\rho_{s}^{2})}{(1 + \rho_{s})(1 + \rho_{s}^{2})} = \dfrac{11}{15}$. \\[0.75em]
The ergodicity condition is not satisfied here, and the queue length will blow up. This is illustrated in the simulation plot of instantaneous queue length below:
\begin{figure}[htp]
	\includegraphics[scale = 0.35]{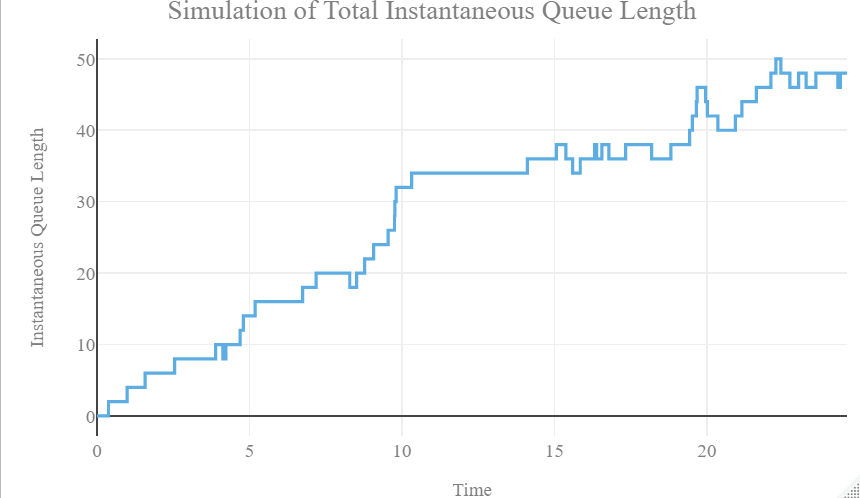}
	\centering
	\caption{Instantaneous queue length for unstable queue}
	\label{fig7}	
\end{figure}
\begin{figure}[htp]
	\includegraphics[scale = 0.35]{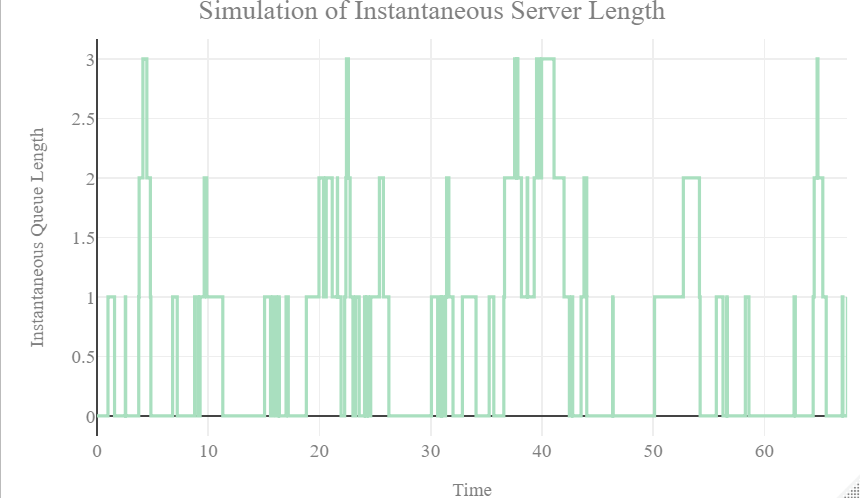}
	\centering
	\caption{Instantaneous server queue length for unstable queue}
	\label{fig8}	
\end{figure}
\begin{figure}[htp]
	\includegraphics[scale = 0.35]{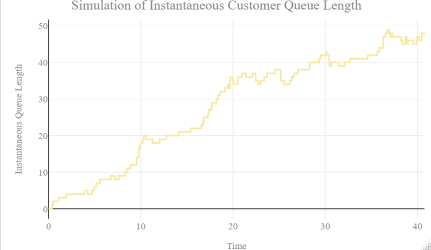}
	\centering
	\caption{Instantaneous customer queue length for unstable queue}
	\label{fig9}	
\end{figure}
\section{Conclusions}\label{sec_5}
In this paper, we present a new kind of queueing system, in which servers arrive at the queue alongside customers with an arrival distribution of their own. An observer to the queue cannot distinguish between a server and a customer, and hence, the queue initially appears to grow bigger. When servers start serving the customers, the size of the queue decreases. We derive equations for this model, find important measures, and demonstrate its behaviour with numerical simulations.

\bigskip
{\bf Acknowledgements.} We acknowledge funding and support  from MITACS Global Internship program.

\end{document}